\documentclass[eqno,12pt]{article}
\usepackage{pstricks}
\usepackage{amsmath}
\usepackage{amsfonts}
\usepackage{amssymb}
\usepackage{graphicx}
\usepackage{epsfig}
\begin{document}
%%%%%%%%%%%%%%%%%%%%%%%%%%%%%%%%%%%%%%%%
\setlength{\leftmargini}{.5\leftmargini}

\newcommand{\ls}[1]
   {\dimen0=\fontdimen6\the\font \lineskip=#1\dimen0
\advance\lineskip.5\fontdimen5\the\font \advance\lineskip-\dimen0
\lineskiplimit=.9
\lineskip \baselineskip=\lineskip
\advance\baselineskip\dimen0 \normallineskip\lineskip
\normallineskiplimit\lineskiplimit \normalbaselineskip\baselineskip
\ignorespaces }

%\ls{1} % single space 
%\ls{2} % double space
%\ls{1.6}
%\ls{1.8}

%TO ADD THE FULL-STOP AFTER THE SECTION NO. E.G. 1. INTRODUCTION

\def\thepart{\Roman{part}} 
\def\thesection {\arabic{section}}
\def\thesubsection {\thesection\arabic{subsection}.}
\def\thesubsubsection {\thesubsection\arabic{subsubsection}.}
\def\theparagraph {\thesubsubsection\arabic{paragraph}.}
\def\thesubparagraph {\theparagraph\arabic{subparagraph}.}

    \setcounter{equation}{0}

\newcommand{\hsp}{{\hspace*{\parindent}}}

%%%%%%%%%%% Equations/theorems/etc, style (1.1) and 1.1
\renewcommand{\theequation}{\thesection.\arabic{equation}}  
\newtheorem{theorem}{Theorem}[section]
\newtheorem{lemma}[theorem]{Lemma}
\newtheorem{proposition}[theorem]{Proposition}
\newtheorem{remark}[theorem]{Remark}
\newtheorem{problem}{Problem}[section]
\newtheorem{definition}{Definition}[section]
\newtheorem{corollary}{Corollary}[section]
\newtheorem{example}{Example}[section]
\newtheorem{conjecture}{Conjecture}[section]
\newtheorem{algorithm}{Algorithm}[section]
\newtheorem{exercise}{Exercise}[section]
%%%%%%%%%%%

% put a period after theorem and theorem-like numbers
\makeatletter
\def\@begintheorem#1#2{\it \trivlist \item[\hskip \labelsep{\bf #1\
#2.}]}
\makeatother

\def\P{{\mathbb P}}
\def\E{{\mathbb E}}

\newcommand{\be}{\begin{equation}}
\newcommand{\ee}{\end{equation}}

\newcommand{\bea}{\begin{eqnarray}}
\newcommand{\eea}{\end{eqnarray}}

\newcommand{\beq}[1]{\begin{equation}\label{#1}}
\newcommand{\eeq}{\end{equation}}
\newcommand{\req}[1]{(\ref{#1})}

\newcommand{\beqn}[1]{\begin{eqnarray}\label{#1}}
\newcommand{\eeqn}{\end{eqnarray}}

\newcommand{\beaa}{\begin{eqnarray*}}
\newcommand{\eeaa}{\end{eqnarray*}}

\newcommand{\eq}[1]{(\ref{#1})}

\def\le{\leq}
\def\ge{\geq}
\def\lt{<}
\def\gt{>}

\newcommand{\lip}{\langle}
\newcommand{\rip}{\rangle}
\newcommand{\uu}{\underline}
\newcommand{\oo}{\overline}
\newcommand{\La}{\Lambda}
\newcommand{\la}{\lambda}
\newcommand{\eps}{\varepsilon}
\newcommand{\vp}{\varphi}

\newcommand{\dint}{\displaystyle\int}
\newcommand{\dsum}{\displaystyle\sum}
\newcommand{\dfr}{\displaystyle\frac}
\newcommand{\bige}{\mbox{\Large\it e}}

%%%%%%%%%%%
\def\NN{{\mathbb N}}
\def\integers{{\mathbb Z}} 
\def\rationals{{\mathbb Q}} 
\def\reals{{\mathbb R}} 
\def\reald{{\reals^d}} 
\def\naturals{{\mathbb N}} 
%%%%%%%%%%%

\newcommand{\degree}{{\scriptscriptstyle \circ }}
\newcommand{\dfn}{\stackrel{\triangle}{=}}
\def\complex{\mathop{\raise .45ex\hbox{${\bf\scriptstyle{|}}$}
     \kern -0.40em {\rm \textstyle{C}}}\nolimits}
\def\hilbert{\mathop{\raise .21ex\hbox{$\bigcirc$}}\kern -1.005em
     {\rm\textstyle{H}}} %Hilbert space

\newcommand{\calA}{{\cal A}}
\newcommand{\calC}{{\cal C}}
\newcommand{\calD}{{\cal D}}
\newcommand{\calF}{{\cal F}}
\newcommand{\calL}{{\cal L}}
\newcommand{\calM}{{\cal M}}
\newcommand{\calP}{{\cal P}}
\newcommand{\calX}{{\cal X}}

\newcommand{\Prob}{{\rm Prob\,}}
\newcommand{\sinc}{{\rm sinc\,}}
\newcommand{\ctg}{{\rm ctg\,}}
\newcommand{\ifff}{\mbox{\ if and only if\ }}
\newcommand{\proof}{\noindent {\bf Proof:\ }}
\newcommand{\remarks}{\noindent {\bf Remarks:\ }}
\newcommand{\note}{\noindent {\bf Note:\ }}

\newcommand{\boldx}{{\bf x}}
\newcommand{\boldX}{{\bf X}}
\newcommand{\boldy}{{\bf y}}
\newcommand{\uux}{\uu{x}}
\newcommand{\uuY}{\uu{Y}}

\newcommand{\liml}{\underline{\lim}_{l\rightarrow \infty}}
\newcommand{\limn}{\underline{\lim}_{n \rightarrow \infty}}
\newcommand{\limN}{\lim_{N \rightarrow \infty}}
\newcommand{\limr}{\lim_{r \rightarrow \infty}}
\newcommand{\limd}{\lim_{\delta \rightarrow \infty}}
\newcommand{\limM}{\lim_{M \rightarrow \infty}}
\newcommand{\limsupn}{\limsup_{n \rightarrow \infty}}
\newcommand{\liminfn}{\liminf_{n \rightarrow \infty}}

\newcommand{\imii}{\int_{-\infty}^{\infty}}
\newcommand{\imix}{\int_{-\infty}^x}
\newcommand{\ioi}{\int_o^\infty}

\newcommand{\ARROW}[1]
  {\begin{array}[t]{c}  \longrightarrow \\[-0.2cm] \textstyle{#1}
\end{array} }

\newcommand{\AR}
 {\begin{array}[t]{c}
  \longrightarrow \\[-0.3cm]
  \scriptstyle {j\rightarrow \infty}
  \end{array}}
\newcommand{\ARn}
 {\begin{array}[t]{c}
  \longrightarrow \\[-0.3cm]
  \scriptstyle {n\rightarrow \infty}
  \end{array}}
\newcommand{\ARm}
 {\begin{array}[t]{c}
  \longrightarrow \\[-0.3cm]
  \scriptstyle {m\rightarrow \infty}
  \end{array}}
\newcommand{\ARN}
 {\begin{array}[t]{c}
  \longrightarrow \\[-0.3cm]
  \scriptstyle {N\rightarrow \infty}
  \end{array}}
\newcommand{\ARk}
 {\begin{array}[t]{c}
  \longrightarrow \\[-0.3cm]
  \scriptstyle {k\rightarrow \infty}
  \end{array}}
\newcommand{\pile}[2]
  {\left( \begin{array}{c}  {#1}\\[-0.2cm] {#2} \end{array} \right) }

\newcommand{\floor}[1]{\left\lfloor #1 \right\rfloor}

%for doing boldface subscripts etc., e.g. $G_{\mmbox{\boldx}}$
\newcommand{\mmbox}[1]{\mbox{\scriptsize{#1}}}

%fraction with round brackets
\newcommand{\ffrac}[2]{\left( \frac{#1}{#2} \right)}
\newcommand{\one}{\frac{1}{n}\:}
\newcommand{\half}{\frac{1}{2}\:}

%qed
\def\squarebox#1{\hbox to #1{\hfill\vbox to #1{\vfill}}}
\newcommand{\qed}{\hspace*{\fill}
\vbox{\hrule\hbox{\vrule\squarebox{.667em}\vrule}\hrule}\smallskip}

%%%%%%%%%%% Renumbering for each section
\newcommand{\eqnsection}{
\renewcommand{\theequation}{\thesection.\arabic{equation}}
    \makeatletter
    \csname  @addtoreset\endcsname{equation}{section}
    \makeatother}
\eqnsection
%%%%%%%%%%% 

%%%%%%%%%%%%%%%%%%%%%%%%%%%%%%%%%%%%%%%%%%%%
%%%%%%%%%%%%%%%%%%%%%%%%%%%%%%%%%%%%%%%%%%%%
%%%%%%%%%%%%%%%%%%%%%%%%%%%%%%%%%%%%%%%%%%%%
%%%
%%%          Text starts here
%%%
%%%%%%%%%%%%%%%%%%%%%%%%%%%%%%%%%%%%%%%%%%%%
%%%%%%%%%%%%%%%%%%%%%%%%%%%%%%%%%%%%%%%%%%%%

\ls{1.5}

\vglue60pt

\begin{centering}
{\Large\bf The infinite valley}

\vskip5pt

{\Large\bf for a recurrent random walk in random environment }

\bigskip
\bigskip
\bigskip

{\sc Nina Gantert}\footnote{CeNos Center for Nonlinear Science, and Institut f\"ur
Mathematische Statistik, Fachbereich Ma\-thema\-tik und Informatik, Universit\"at M\"unster, Einsteinstr.\ 62, D-48149 M\"unster, Germany, email:
gantert@math.uni-muenster.de  Research partially
supported by the European program RDSES}, {\sc Yuval Peres}\footnote{Department of Statistics, University of California at Berkeley, Berkeley, CA 94720, email: peres@stat.berkeley.edu} and  
{\sc Zhan
Shi}\footnote{Laboratoire de Probabilit\'es et Mod\`eles Al\'eatoires, Universit\'e
Paris VI, 4 place Jussieu, F-75252 Paris Cedex 05, France, email:
zhan.shi@upmc.fr}

{\tt This version: February 25, 2009}

\end{centering}

\bigskip
\bigskip
\bigskip

\noindent
{\sc Abstract:} We consider a one-dimensional recurrent random walk in random environment (RWRE). We show that the -- suitably centered -- empirical distributions of the RWRE converge weakly to a certain limit law which describes the stationary distribution of a random walk in an infinite valley. The construction of the infinite valley goes back to Golosov, see \cite{golosov}. As a consequence, we show weak convergence for both the maximal local time and the self-intersection local time of the RWRE and also determine 
the exact constant in the almost sure upper limit of the maximal local 
time.

\noindent
{\sc R\'esum\'e:}
Nous prouvons que les mesures empiriques d'une marche
al\'eatoire unidimensionnelle en environnement al\'eatoire
convergent \'etroitement vers la loi stationnaire d'une marche
al\'eatoire dans une vall\'ee infinie. La construction de cette
vall\'ee infinie revient \`a Golosov, voir \cite{golosov}. En
applications, nous obtenons la convergence \'etroite du maximum
des temps locaux et du temps local d'intersections de la marche
al\'eatoire en environnement al\'eatoire ; de plus, nous
identifions la constante repr\'esentant la ``limsup" presque
s\^ure du maximum des temps locaux.


\bigskip


\bigskip

\ls{1.5}

\section{Introduction and statement of the results}
\label{s:intro}

Let $\omega =(\omega_x)_{x \in \integers}$ be a collection of i.i.d. random variables taking values in 
$(0,1)$ and let $P$ be the distribution of $\omega$. For each $\omega \in \Omega= (0,1)^{\integers}$, we define the random walk in random environment (abbreviated RWRE) as the time-homogeneous Markov chain $(X_n)$ taking values in $\integers_+$, with transition probabilities $P_\omega(X_{n+1} = 1|X_n = 0)=1$, $P_\omega [X_{n+1} = x+1 \, |\, X_n =x] = \omega_x = 1-P_\omega [X_{n+1} = x-1 \, | \, X_n =x]$ for $x >0$, and $X_0 =0$. We equip $\Omega$ with its Borel $\sigma$-field
$\cal{F}$ and $\integers^\naturals$ with its Borel $\sigma$-field $\cal{G}$. The distribution of $(\omega, (X_n))$ is the probability measure $\P$ on $\Omega \times \integers^\naturals$ defined by $\P[F\times G] = \int\limits_F P_\omega[G] P(d\omega)$, $F \in \cal{F}$, $G\in\cal{G}$. Let $\rho_i= \rho_i(\omega): = (1-\omega_i)/\omega_i$. We will always assume that 
\begin{eqnarray}
    \int \log \rho_0(\omega) \, P(d\omega) =0\, ,
\label{rec}
\end{eqnarray}    
\begin{eqnarray}
    P[\delta \leq \omega_0  \leq 1-\delta] = 1
 \hbox{ for some } \delta \in (0,1),
    \label{bddsupp} 
\end{eqnarray}
\begin{eqnarray}
\hbox{\rm Var} (\log \rho_0) > 0\, .
    \label{det}
\end{eqnarray}
The first assumption, as shown in \cite{solomon}, implies that for $P$-almost all $\omega$, the Markov chain $(X_n)$ is recurrent, the second is a technical assumption which could probably be relaxed but is used in several places, and the third assumption excludes the deterministic case.
\noindent 
Usually, one defines in a similar way the RWRE on the integer axis, 
but for simplicity, we stick to the RWRE on the positive integers; see Section \ref{s:remarks} for the conjectured results for the usual RWRE model.
A key property of recurrent RWRE is
its strong localization: under the assumptions above, Sinai \cite{Sinai} showed 
that $X_n/(\log n)^2$  converges in distribution.
A lot more is known about this model; we refer to the  survey 
by Zeitouni~\cite{zeitouni} for 
limit theorems, large deviations results, and for further references. 

 Let $\xi(n,x): = |\{0 \le j \le n: X_j =x\}|$
denote the local time of the RWRE in $x$ at time $n$ and  $\xi^*(n) :=
\sup_{x\in \integers} \xi(n,x)$ the maximal local time at time $n$. It was
shown in \cite{r90} and \cite{zhan1} that
\begin{equation}\label{isconst} 
\limsup_{n\to \infty} \, \frac{\xi^*(n)}{ n} >0
    \qquad \hbox{\rm $\P$-a.s.}
\end{equation}
(Clearly this $\limsup$ is at most $1/2$.)
%new
\noindent In addition, a 0--1 law (see \cite{mama}) says that $\limsup_{n \to \infty}\frac{\xi^*(n)}{n}$ is $\P$-almost surely a constant. The constant however was not known. We give its value in Theorem \ref{t:LIL}.

Our main result (Theorem \ref{empdistr} below) shows weak convergence for the process $(\xi(n,x), \, x\in\integers)$ -- after a suitable normalization -- in a function space. In particular, it will imply the following theorem.

\medskip
\begin{theorem}
 \label{t:LIL}
Let $M:= \sup \{ s: \, s\in \hbox{\rm 
 supp}(\omega_0)\} \in (\frac{1}{2}, \, 1]$ and
$w:= \inf \{ s: \, s\in \hbox{\rm 
 supp}(\omega_0)\} \in [0, \, \frac{1}{2})$.
Then 
\begin{equation}\label{theconst}
\limsup_{n\to \infty} \, \frac{\xi^*(n)}{n} =\frac{(2M-1)(1-2w)}{2(M-w) \min\{ M, 1-w\}}, \qquad \P\hbox{\rm -a.s.}
\end{equation}
In particular, if $M= 1-w$, we have
\begin{equation}
     \limsup_{n\to \infty} \, \frac{\xi^*(n)}{n} =  
     \frac{2M-1}{2M}, \qquad \P\hbox{\rm -a.s.}
 \end{equation}

\end{theorem}

Define the potential $V=(V(x), \, x\in\integers)$ by
$$
V(x): = \begin{cases}
\sum\limits_{i=1}^x \log \rho_i , &  x > 0 \\
0, & x = 0\\ 
-\sum\limits_{i=x+1}^{0}\log\rho_{i}, & x < 0
\end{cases}
$$

\noindent and $C_{(x, x+1)}:= \exp(-V(x))$. For each $\omega$, the Markov chain is an electrical network in the sense of 
\cite{doylesnell}, where $C_{(x, x+1)}$ is the conductance of the bond $(x, x+1)$. In particular, $\mu(x) : = \exp(-V(x-1)) + \exp(-V(x)), \, x > 0$, $\mu(0) =1$ is a (reversible) invariant measure for the Markov chain.

Let ${\widetilde V}= ({\widetilde V}(x), \, x\in \integers)$ be a collection of random variables distributed as $V$ conditioned to stay non-negative for $x > 0$ and strictly positive for $x<0$.
Due to (\ref{det}), such a distribution is well-defined, see for example Bertoin~\cite{bertoin} or Golosov~\cite{golosov}.
Moreover, it has been shown that
\begin{equation}\label{expsumm}
\sum_{x\in \integers} \exp(-{\widetilde V}(x)) < \infty,
\end{equation}
see \cite{golosov}, p. 494. For each realization of $({\widetilde V}(x), \, x\in \integers)$ consider the corresponding Markov chain on $\integers$, which is an electrical network with conductances $\widetilde{C}_{(x, x+1)}:= \exp(-{\widetilde V}(x))$. Intuitively, this Markov chain is a random walk in the
 ``infinite valley" given by ${\widetilde V}$. As usual, 
${\widetilde \mu}(x) : = \exp(-{\widetilde V}(x-1)) + \exp(-{\widetilde V}(x)), \, x \in \integers $ is a reversible measure for this Markov chain. But, due to (\ref{expsumm}), we can normalize ${\widetilde \mu}$ to get a reversible probability measure $\nu$, defined by
\begin{equation}
\label{nudef}
    \nu(x) 
 := \frac{\exp(- {\widetilde V} (x-1))+ \exp( -
    {\widetilde V} (x))}{
    2\sum\limits_{x \in \integers} \exp(-
    {\widetilde V} (x))}\, ,\qquad x \in \integers.
\end{equation}

\noindent Note that in contrast to the original RWRE which is null-recurrent, the random walk in the ``infinite valley" given by ${\widetilde V}$ is positive recurrent.

Let $\ell^1$ be the space of real-valued sequences $\ell := \{\ell(x), \, x\in \integers\}$ satisfying $\| \ell \| := \sum_{x\in \integers} \vert \ell(x)\vert <\infty$.
Let 
$$
c_n: = \min\left\{x\geq 0: V(x) - \min\limits_{0\leq y \leq x} V(y)\geq \log n + (\log n)^{1/2}\right\}
$$
and 
$$
b_n : = \min\left\{x\geq 0: V(x) = \min\limits_{0\leq y \leq c_n} V(y)\right\}\, .
$$
\vfill\break
\begin{figure}[h]
 \begin{pspicture}(-1,-5)(15,5)
  \psline{->}(-1,0)(14,0)
  \psline{->}(0,-4)(0,5)

  \psline{*-*}(0,0)(0.4,0.5)
  \psline{*-*}(0.4,0.5)(0.8,-0.5)
  \psline{*-*}(0.8,-0.5)(1.2,1)
  \psline{*-*}(1.2,1)(1.6,1.5)
  \psline{*-*}(1.6,1.5)(2,-0.2)
  \psline{*-*}(2,-0.2)(2.4,-1)
  \psline{*-*}(2.4,-1)(2.8,-0.5)
  \psline{*-*}(2.8,-0.5)(3.2,-1.7)
  \psline{*-*}(3.2,-1.7)(3.6,-1.5)
  \psline{*-*}(3.6,-1.5)(4,-2)
  \psline{*-*}(4,-2)(4.4,-3.5)
  \psline{*-*}(4.4,-3.5)(4.8,-3.1)
  \psline{*-*}(4.8,-3.1)(5.2,-3.8)
  \psline{*-*}(5.2,-3.8)(5.6,-2.8)
  \psline{*-*}(5.6,-2.8)(6,-3.1)
  \psline{*-*}(6,-3.1)(6.4,-2.9)
  \psline{*-*}(6.4,-2.9)(6.8,-3.8)
  \psline{*-*}(6.8,-3.8)(7.2,-3.5)
  \psline{*-*}(7.2,-3.5)(7.6,-2.4)
  \psline{*-*}(7.6,-2.4)(8,-3.8)
  \psline{*-*}(8,-3.8)(8.4,-2.2)
  \psline{*-*}(8.4,-2.2)(8.8,-1.2)
  \psline{*-*}(8.8,-1.2)(9.2,-1.6)
  \psline{*-*}(9.2,-1.6)(9.6,0.4)
  \psline{*-*}(9.6,0.4)(10,-0.4)
  \psline{*-*}(10,-0.4)(10.8,2)
  \psline{*-*}(10.8,2)(11.2,1.7)
  \psline{*-*}(11.2,1.7)(11.6,3.4)
  \psline{*-*}(11.6,3.4)(12,4.4)
  \psline[linestyle=dotted]{-}(11,4)(13,4)
  \psline[linestyle=dotted]{-}(12,4.4)(12,0)
  \psline[linestyle=dotted]{-}(4.2,-3.8)(13,-3.8)
  \psline[linestyle=dotted]{-}(5.2,-3.8)(5.2,0)
  \psline[linestyle=dotted]{<-}(12.5,4)(12.5,-1.2)
  \psline[linestyle=dotted]{<-}(12.5,-3.8)(12.5,-1.8)

  \rput(5.2,0.3){$b_n$}
  \rput(12,-0.3){$c_n$}
  \rput(-0.3,0.3){0}
  \rput(0.5,5){$V(x)$}
  \rput(14,-0.3){$x$}
  \rput(12.5,-1.5){$\log n + (\log n)^{1/2}$}

  \rput(6, -5){Figure 1: definition of $b_n$}
  
 \end{pspicture}
\end{figure}

\medskip

We will consider  $\{ \xi(n, x) , \, x \in \integers \}$ as a random element of $\ell^1$ (of course, $ \xi(n, y)= 0 $ for $y < 0$).
Here is our main result.
\medskip
\begin{theorem}
 \label{empdistr}
 Consider $\{ \xi(n, x) , \, x \in \integers \}$ as 
 a random element of $\ell^1$ under the probability 
 $\P$. Then, \begin{equation}\label{mainconv}
     \left\{\frac{\xi(n, b_n+x )}{n}\; , x \in   
     \integers \right\} \quad {\buildrel law \over 
     \longrightarrow} \quad \nu \, , \qquad n\to 
     \infty,
 \end{equation}
 where ${\buildrel law \over \longrightarrow}$ 
 denotes convergence in distribution. In other words, the distributions of \\
$\left\{\frac{\xi(n, b_n+x )}{n}\; , x \in   \integers \right\}$ converge weakly to the distribution of $\nu$ (as probability measures on $\ell^1$).
\end{theorem}

The following corollary is immediate.
\medskip

\begin{corollary}
 \label{funct}
 For each continuous functional $f: \ell^1 \to  
 \reals$ which is shift-invariant, we have
 \begin{equation}\label{fucoll}
     f\left(\left \{\frac{\xi(n, x)}{n}\, , x  
     \in \integers\right \}\right) \quad {\buildrel law 
     \over \longrightarrow} \quad f\left(\{\nu(x)\, , x 
     \in \integers\}\right)\, .
 \end{equation}
\end{corollary}

\medskip

\noindent {\bf Example 1 (Self-intersection local time).} Let 
$$
f\left(\{\ell(x), x \in \integers\}\right) = \sum\limits_{x \in \integers} \ell(x)^2\, .
$$

\noindent Corollary \ref{funct} yields that
\begin{equation}\label{sil}
\frac{1}{n^2}\sum\limits_{x \in \integers}\xi(n, x)^2 \quad {\buildrel law \over \longrightarrow} \quad \sum\limits_{x \in \integers} \nu(x)^2 \, .
\end{equation}
This confirms Conjecture 7.4 in \cite{zhan2}, at least for the RWRE on the positive integers; for the RWRE on the integer axis, see 
Section \ref{s:remarks}.
\bigskip

\noindent {\bf Example 2 (Maximal local time).} Let 
$$
f\left(\{\ell(x), x \in \integers\}\right) = \sup\limits_{x \in \integers} \ell(x)\, .
$$

\noindent Corollary \ref{funct} yields that
\begin{equation}\label{maxloc}
\frac{\xi^*(n)}{n} \quad {\buildrel law \over \longrightarrow} \quad \sup\limits_{x \in \integers} \nu(x)\, .
\end{equation}

\section{Proof of Theorem \ref{empdistr}}
\label{s:proof}

The key lemma is

\medskip
\begin{lemma}
 \label{l:ficonf}
 For any $K \in \NN$,
 \begin{equation}
 \left\{ \frac{\xi(n, \, b_n+x)}{n}, -K \leq x \leq K\right\} {\buildrel law \over \longrightarrow}
\left\{\nu(x), -K \leq x \leq K\right\}.
    \label{ficonf}
\end{equation}
\end{lemma}

Given Lemma \ref{l:ficonf}, the proof of Theorem \ref{empdistr} is straightforward. It suffices to 
show that for every function 
$f:\ell^1 \to \reals$, $f$  bounded and uniformly continuous, we have
\begin{equation}
\label{want}
\E\left[f\left(\left\{\frac{\xi(n, \, b_n+x)}{n}, x \in \integers\right\}\right)\right] \to 
E\left[f\left(\left\{\nu(x), x \in \integers\right\}\right)\right], \quad n \to \infty \, .
\end{equation}
For $\ell \in \ell^1$, define $\ell_K$ by $\ell_K(x) = \ell(x)I_{\{-K \leq x \leq K\}}$  ($x \in \integers $). For 
$f:\ell^1 \to \reals$, define $f_K$ by $f_K(\ell) = f(\ell_K)$ ($\ell \in \ell^1$).
Due to Lemma \ref{ficonf},
\begin{equation}
\label{fMconf}
\E\left[f_K\left(\left\{\frac{\xi(n, \, b_n+x)}{n}, x \in \integers\right\}\right)\right] \to 
E\left[f_K\left(\left\{\nu(x), x \in \integers\right\}\right)\right], \quad n \to \infty\, .
\end{equation}
Fix $\varepsilon > 0$. Since $f$ is uniformly continuous, there is ${\widetilde \delta}$ such that we have for $\ell \in \ell^1$
\begin{equation}
\|\ell_K - \ell\|_1 \leq {\widetilde \delta} \Rightarrow \vert f_K(\ell) -f(\ell)\vert \leq \varepsilon .
\end{equation}
Hence
\begin{equation}
\E\left[\left\vert f_K \left(\left\{\frac{\xi(n, \, b_n+x)}{n}, x \in \integers\right\}\right) - f\left(\left\{\frac{\xi(n, \, b_n+x)}{n}, x \in \integers\right\}\right)\right\vert\right]\leq \varepsilon
\end{equation}
provided that
\begin{equation}
\label{fest}
\E\left[\sum\limits_{x: |x| > K}\frac{\xi(n,b_n+x )}{n}\right]\, \leq {\widetilde \delta}\, .
\end{equation}
Due to Lemma \ref{l:ficonf}, $\E\left[\sum\limits_{x: |x| > K}\frac{\xi(n, \, b_n+x)}{n}\right]$ converges to $E\left[\sum\limits_{x: |x| > K}\nu(x)\right]$. But, for $K$ large enough, $E\left[\sum\limits_{x: |x| > K}\nu(x)\right] \leq {\widetilde \delta}$ since 
$E\left[\nu(\cdot)\right]$ is a probability measure on $\integers$.\ Together with (\ref{fest}) and (\ref{fMconf}), this implies (\ref{want}).
\qed

\noindent {\it Proof of Lemma \ref{l:ficonf}.}\\
We proceed in several steps.\\
(i) Define
\begin{equation}
\label{Tdef}
T(y) := \inf\left\{ n\geq 1: \; X_n = y \right\},
\end{equation}
the first hitting time of $y$. For $\varepsilon > 0$, let $A_n$ denote the event $\{T(b_n) \leq \varepsilon n\}$. Then, 
$\P[A_n] \to 1$ for $n \to \infty$.\\
Proof: see Golosov \cite{golosov}, Lemma 1. \\
(ii) Let $B_n$ denote the event $\{\hbox{the RWRE exits}\,  [0, c_n) \hbox { before time } n\} = \{T(c_n) < n\}$. Then, for $P$-a.a. $\omega$,
$P_\omega[B_n | X_0 = b_n] \to 0$ for $ n \to \infty$.\\
Proof:  Due to \cite{golosov}, Lemma 7, we have for all $m$
$$
P_\omega[T(c_n) < m | X_0 = b_n] \leq {\rm const}\cdot m e^{-\log n - (\log n)^{1/2}} \, .
$$
Taking $m= n$, we conclude that $P_\omega[B_n | X_0 = b_n] \to 0$.\\
(iii) Due to (i) and (ii), we can consider, instead of $P_\omega$, a finite Markov chain $\widetilde P_\omega = \widetilde P_\omega^{(n)}$ started from $ b_n$, in the valley $[0, c_n]$. More precisely, the Markov chain $\widetilde P_\omega$ is the original Markov chain $P_\omega$ with reflection at $c_n$, i.e. for $0< x < c_n$,
$\widetilde P_\omega [X_{n+1} = x+1 \, |\, X_n =x] = \omega_x = 1-\widetilde P_\omega [X_{n+1} = x-1 \, | \, X_n =x]$, $\widetilde P_\omega [X_{n+1} = 1 \, |\, X_n =0] = \widetilde P_\omega [X_{n+1} = c_n -1 \, |\, X_n =c_n] = 1$ and $X_0 =b_n$.
The invariant probability measure $\mu_\omega = \mu_\omega^{(n)}$ of $\widetilde P_\omega$ is given by
$$
\mu_\omega(x): = \begin{cases}
\frac{1}{Z_\omega}\left(e^{-V(x)} +e^{-V(x-1)}\right)  , &  0 < x < c_n \\
\frac{1}{Z_\omega}e^{-V(0)}, & x = 0\\ 
\frac{1}{Z_\omega}e^{-V(c_n-1)} ,& x = c_n
\end{cases}
$$
where $Z_\omega = 2\sum\limits_{x=0}^{c_n-1}e^{-V(x)}$.\\
(iv) Recall (\ref{Tdef}).
We note for further reference that for $0 \leq b < y < i$,
\begin{equation}
\label{hitbebeforei}
    P_\omega\left[T(b) < T(i) | X_0 = y\right] = 
    \sum\limits_{j=y}^{i-1} e^{V(j)}
    \Big(\sum\limits_{j=b}^{i-1} e^{V(j)}\Big)^{-1} .
\end{equation}
This follows from direct computation, using $C_{(x, x+1)} = e^{-V(x)}$, 
see also \cite{zeitouni}, formula (2.1.4).
We now decompose the paths of the Markov chain $\widetilde P_\omega$ into excursions from $b_n$ to $b_n$.
Let $x \in [-b_n, c_n - b_n]$, $x \neq 0$ and  denote by $Y_{b_n,x}$ the
number of visits to $b_n + x$ before returning to $b_n$.
The distribution of $Y_{b_n,x}$ is \lq\lq almost geometric\rq\rq : we have
$$
\widetilde P_\omega[Y_{b_n,x} =m] =
\begin{cases} \alpha(1-\beta)^{m-1}\beta &m=1,2,3, \ldots, \cr
        1-\alpha, &m=0, \cr
\end{cases}
$$
where $\alpha = \alpha_{b_n,x} = \widetilde P_\omega[T(b_n+x) < T(b_n) \, | \, X_0 = b_n]$,
$\beta = \beta_{b_n,x} = \widetilde P_\omega[T(b_n) < T(b_n+x) | X_0 = b_n+x]$. In
particular,
\begin{equation}
\label{expu}
    \widetilde E_\omega[Y_{b_n,x}] = \frac{\alpha}{\beta} = \frac{\mu_\omega(b_n+x)}{\mu_\omega(b_n)}
\end{equation}
\noindent where $\mu_\omega$ is the reversible measure for the Markov chain,
see above. Note that (\ref{expu}) also applies for $x=0$, with $Y_{b_n,0} = 1$.
Now, with $\hbox{Var}_\omega(Y_{b_n,x})$ denoting the variance of $Y_{b_n,x}$ w.r.t. $\widetilde P_\omega$,
\begin{equation}
\label{varest}
\hbox{Var}_\omega(Y_{b_n,x}) = \frac{\alpha (2- \beta - \alpha)}
{\beta^2}
\leq \frac{2}{\beta} \frac{\mu_\omega(b_n+x)}{\mu_\omega(b_n)} \leq \frac{4}{\beta}\, .
\end{equation}
For $x >1$, 
$$
\beta = (1-\omega_{b_n+x})\widetilde P_\omega[T(b_n)< T(b_n+x)| X_0 = b_n +x-1]\, .
$$
Taking into account (\ref{hitbebeforei}) yields
\begin {equation}
\label{betaval}
\beta = (1-\omega_{b_n+x})\Big(\sum\limits_{j=b_n}^{b_n+x-1} e^{V(j) - V(b_n+ x-1)}\Big)^{-1} ,
\end{equation}
and (\ref{betaval}) applies also 
to
$x=1$. In particular, recalling (\ref{bddsupp}), 
there is a constant $a=a(K,\delta)>0$ such that for $1 \leq x \leq K$, uniformly in $n$,
\begin{equation}
    \hbox{ Var}_\omega(Y_{b_n,x}) 
 \leq  a(K, \delta) \, .  
    \label{varvis}
\end{equation}
 In the same way, one obtains (\ref{varvis}) for $-K \leq x \leq 0$.
We note for further reference that due to (\ref{varest}) and (\ref{betaval}), there is a constant $\widetilde{a} = \widetilde{a}(\delta) > 0$ such that for $x \in [-b_n, c_n - b_n]$,
\begin{equation}
    \hbox{ Var}_\omega(Y_{b_n,x}) 
 \leq  \widetilde{a}(\delta)c_n (\log n + (\log n)^{1/2}) \, .  
    \label{varvis2}
\end{equation}
(v) Denote by $k_n$ the number of excursions from $b_n$ to $b_n$ before time $n$. It follows from (\ref{expu}) that the average 
length $\gamma_n$ of an excursion from $b_n$ to $b_n$ under $\widetilde P_\omega$ is given by
$$
\gamma_n = \sum\limits_{y=0}^{c_n} \frac{\mu_\omega(y)}{\mu_\omega(b_n)}\geq 2 \, .
$$
Fix $\varepsilon \in (0,1)$. We show that for $P$-a.a. $\omega$,
\begin{equation}
\label{numbex}
\widetilde P_\omega\left[\left\vert \frac{k_n}{n} - \frac{1}{\gamma_n}\right\vert \geq \varepsilon\right] \to 0, 
\quad n \to \infty\, .
\end{equation}
Proof: Let $Y_{b_n, x}^{(1)}, Y_{b_n, x}^{(2)}, Y_{b_n, x}^{(3)}, \ldots $ be i.i.d. copies of $Y_{b_n, x}$, and denote by $E^{(i)}_n = \sum\limits_{x=-b_n}^{c_n-b_n}Y_{b_n, x}^{(i)}$ the length of the $i$-th excursion from $b_n$ to $b_n$. Then, 
$\{k_n \geq  m\} \subseteq \{\sum\limits_{i=1}^m E^{(i)}_n \leq n\}$ and $\{k_n \leq  m\} \subseteq \{\sum\limits_{i=1}^m E^{(i)}_n \geq n\}$.
Hence,
\begin{equation}
\label{umdreh}
\widetilde P_\omega\left[\left\vert \frac{k_n}{n} - \frac{1}{\gamma_n}\right\vert \geq \varepsilon\right] = 
\widetilde P_\omega\left[\sum\limits_{i=1}^{n(\frac{1}{\gamma_n} + \varepsilon)} E^{(i)}_n \leq n \right] +
\widetilde P_\omega\left[\sum\limits_{i=1}^{n(\frac{1}{\gamma_n} - \varepsilon)} E^{(i)}_n \geq n \right]\, .
\end{equation}
To handle the first term in (\ref{umdreh}), recall $E^{(1)}_n, E^{(2)}_n, \ldots $ are i.i.d. with expectation $\gamma_n$ under $\widetilde P_\omega$ and apply Chebyshev's inequality:
\begin{equation}
\label{cheba}
\widetilde P_\omega\left[\frac{1}{n(\frac{1}{\gamma_n} + \varepsilon)}\sum\limits_{i=1}^{n(\frac{1}{\gamma_n} + \varepsilon)} \left(E^{(i)}_n - \gamma_n\right) \leq  - \frac{\varepsilon\gamma_n^2}{1 + \gamma_n \varepsilon}\right]
\leq \frac{{\rm Var}_\omega(E^{(1)}_n)} {n(\frac{1}{\gamma_n}+ \varepsilon)} \frac{(1 + \gamma_n \varepsilon)^2}{\varepsilon^2\gamma_n^4} \, . 
\end{equation}        
Now, use the inequality ${\rm Var}(\sum\limits_{j=1}^N X_j) \leq N \sum\limits_{j=1}^N {\rm Var}(X_j)$ and (\ref{varvis2}) to get from (\ref{cheba}) that
$$
\widetilde P_\omega\left[\sum\limits_{i=1}^{n(\frac{1}{\gamma_n} + \varepsilon)} E^{(i)}_n \leq  n\right]
\leq \frac{c_n^4 \widetilde{a}(\delta)^2(\log n)^3}{n}
\frac{(1 + \gamma_n \varepsilon)}{\varepsilon^2\gamma_n^3}\leq \frac{c_n^4 \widetilde{a}(\delta)^2(\log n)^3}{n}
\frac{2}{\varepsilon^2}  \, .
$$
Due to Chung's law of the iterated logarithm, \\
$\liminf\limits_{ n\to \infty}\left(\frac{\log\log n}{n}\right)^{1/2}\left(\max_{0 \leq x \leq n} V(x) -  \min_{0 \leq x \leq n} V(x)\right)$ is a.s. a strictly positive constant, and this implies
$c_n^4/n  \to 0$ for $P$-a.a. $\omega$. We conclude that the first term in (\ref{umdreh}) goes to $0$ for a.a. $\omega$. The second term in (\ref{umdreh}) is treated in the same way.\\
Further, due to Chebyshev's inequality, (\ref{expu}) and (\ref{varvis}), we have for $-K \leq x \leq K$
\begin{equation}
\label{numbvis}
\widetilde P_\omega\left[\left\vert \frac{1}{k_n}\sum\limits_{i=1}^{k_n} Y_{b_n, x}^{(i)} - 
\frac{\mu_\omega(b_n+x)}{\mu_\omega(b_n)}\right\vert \geq  \varepsilon\right]\to 0, \quad n \to \infty\, .
\end{equation}
(vi) Define, for $-K \leq x \leq K$, $\rho_{n, \omega}$ by
\begin{equation}
\rho_{n, \omega}(x) = \frac{\xi(n+ T(b_n), \,b_n+x)}{n} - \frac{\xi(T(b_n), \,b_n+x)}{n}\quad, -K \leq x \leq K\, .
\end{equation}
For $x=0$, we have $n\rho_{n, \omega}(x) = k_n$.\\
We now estimate $\rho_{n, \omega}(x)$ for $-K \leq x \leq K$, $x \neq 0$:
\begin{equation}
\frac{k_n}{n}\frac{1}{k_n}\sum\limits_{i=1}^{k_n} Y_{b_n, x}^{(i)} \leq \rho_{n, \omega}(x)
\leq \frac{k_n+1}{n}\frac{1}{k_n+1}\sum\limits_{i=1}^{k_n+1} Y_{b_n, x}^{(i)} \, .
\end{equation}
We conclude from (\ref{numbex}) and (\ref{numbvis}) that for $-K \leq x \leq K$,
\begin{equation}
\widetilde P_\omega\left[\left\vert \rho_{n, \omega}(x) - \frac{1}{\gamma_n}\frac{\mu_\omega(b_n+x)}
{\mu_\omega(b_n)}\right\vert\geq \varepsilon \right] \to 0, \quad n \to \infty\, .
\end{equation}
(vii) We show that $\left\{\frac{1}{\gamma_n}\frac{\mu_\omega(b_n+x)}{\mu_\omega(b_n)}, -K \leq x \leq K\right\}$ converges in distribution to\\
 $\left\{\nu(x), -K \leq x \leq K\right\}$.\\
Proof: Due to \cite{golosov}, Lemma 4, the finite-dimensional distributions of $\left\{V( b_n + x)- V(b_n)\right\}_{x \in \integers}$ converge to those of $\left\{\widetilde V( x)\right\}_{x \in \integers}$. Therefore, setting $\gamma_n^{(N)}: = \sum\limits_{x=b_n-N}^{b_n+N} \frac{\mu_\omega(x)}{\mu_\omega(b_n)}$, we have for each $N$ that 
$$
\left\{\frac{1}{\gamma_n^{(N)}} \frac{\mu_\omega(b_n+x)}{\mu_\omega(b_n)}, -K \leq x \leq K\right\}
$$
converges in distribution to 
$$
\left\{\frac{e^{-\widetilde V(x)} +e^{-\widetilde V(x-1)}}{\sum\limits_{x= -N}^N e^{-\widetilde V(x)} +e^{-\widetilde V(x-1)}} , -K \leq x \leq K\right\}\, .
$$
It remains to show that for $N$ large enough, $\gamma_n^{(N)}$ and $\gamma_n$ are close for $n \to \infty$ and that, for $N$ large enough,
$\sum\limits_{x= -N}^N \left(e^{-\widetilde V(x)} +e^{-\widetilde V(x-1)}\right)$ is close to $\sum\limits_{x= -\infty}^\infty \left(e^{-\widetilde V(x)} +e^{-\widetilde V(x-1)}\right)$. Note that $\gamma_n^{(N)} \leq \gamma_n$ and
\begin{equation}
\label{smala}
P\left[\frac{\gamma_n^{(N)}}{\gamma_n} \geq 1-\varepsilon\right] \to P\left[\frac{\sum\limits_{x= -N}^N \left(e^{-\widetilde V(x)} +e^{-\widetilde V(x-1)}\right)}{\sum\limits_{x= -\infty}^\infty \left(e^{-\widetilde V(x)} +e^{-\widetilde V(x-1)}\right)} \geq 1-\varepsilon\right]
\end{equation}
for $n \to \infty$. But, due to (\ref{expsumm}), $\left\{ \sum\limits_{x= -N}^N \Bigl(e^{-\widetilde V(x)} +e^{-\widetilde V(x-1)}\Bigr)\geq(1-\varepsilon) \sum\limits_{x= -\infty}^\infty \Bigl(e^{-\widetilde V(x)} +e^{-\widetilde V(x-1)}\Bigr)\right\}$ increases to the whole space
and we conclude that the last term in (\ref{smala}) goes to $1$ for $N \to \infty$.\\

Putting together (i) - (vii), we arrive at (\ref{ficonf}).\qed

\section{Proof of Theorem \ref{t:LIL} }
\label{s:constant}

We know now that the laws of $\frac{\xi^*(n)}{ n}$ converge to the law of $\nu^*$ (see (\ref{maxloc})), where
\begin{equation}\label{stardef}
\nu^* : = \sup_{x\in \integers}{\nu}(x) =  \sup_{x\in \integers} \frac{\exp(- {\widetilde V} (x-1))+ \exp( - {\widetilde V} (x))}{ 2\sum\limits_{x \in \integers} \exp(- {\widetilde V} (x))}  \, .
\end{equation}
We will show that 
\begin{equation}\label{coco}
     \limsup_{n\to \infty} \, \frac{\xi^*(n)}{ n} =  c, \qquad \P\hbox{\rm -a.s.},
 \end{equation}
where
\begin{equation}\label{cocon}
c = \sup \{z: z \in \hbox{\rm supp}\,  \nu^* \}\, .
\end{equation}
In fact, we get
\begin{equation}\label{susu}
     \limsup_{n\to \infty} \, \frac{\xi^*(n)}{n} \geq c
\qquad \P\hbox{\rm -a.s.}
 \end{equation}
from the following general
lemma, whose proof is easy (and therefore omitted).
\begin{lemma}
\label{easyboundsupp}
Let $(Y_n), n=1,2, \ldots $ be a sequence of random variables such that $(Y_n)$ converges in distribution to a 
random variable $Y$ and $ \limsup_{n\to \infty} Y_n$ is constant almost surely.
Then, 
\begin{equation}
\limsup_{n\to \infty} Y_n \geq \sup \{z: z \in \hbox{\rm supp}\, Y \}\, .
\end{equation}
\end{lemma}
To show that
\begin{equation}
     \limsup_{n\to \infty} \, \frac{\xi^*(n)}{n} \leq c 
\qquad \P\hbox{\rm -a.s.},
 \end{equation}
we will use a coupling argument with an environment $\overline \omega$ whose potential $\overline V$ will be shown to achieve the supremum in
(\ref{cocon}).
Define the environment $\overline \omega$ as follows.
$$
{\overline \omega}_x : = \begin{cases}
 w, &  x > 0 \\
M, & x \leq 0.\\
\end{cases}
$$
For the corresponding potential $\overline V$, $\exp(- {\overline V(x)})$ is given as follows: 
$$
\exp(- {\overline V(x)})  = \begin{cases}
\left(\frac{w}{1-w}\right)^x,  &  x > 0 \\
1, & x = 0\\ 
\left(\frac{M}{1-M}\right)^x,   & x < 0. 
\end{cases}
$$
For any fixed environment $\omega \in [w, M]^\integers$, $Q_\omega$ defines a Markov chain on the integers with transition probabilities $Q_\omega [X_{n+1} = x+1 \, |\, X_n =x] = \omega_x = 1-Q_\omega [X_{n+1} = x-1 \, | \, X_n =x]$, $x \in \integers$, and $X_0 =0$.
Then we have the following lemma.

\begin{lemma}
\label{coupling}
Let the environment $\overline{\omega}$ be defined as above. Assume $M \leq 1-w$. For {\sl all} $\omega \in [w, M]^\integers$, all $x \in \integers$ and all $n \in \naturals$, the distribution of $\xi(n,x)$ with respect to $Q_\omega$ is stochastically dominated by the distribution of $\xi(n,0)$ with respect to $Q_{\overline{\omega}}$. In particular, 
\begin{equation}\label{nubar}
\sup \{z: z \in \hbox{\rm supp} \, \sup_{x\in \integers} {\nu}(x)\}
= {\overline\nu}(0):=
\frac{\exp(- {\overline V} (-1))+ \exp( - {\overline V} (0))}{ 2\sum\limits_{x \in \integers} \exp(- {\overline V} (x))} . 
\end{equation}
\end{lemma}

\noindent {\it Proof of Lemma \ref{coupling}.} First step: we show that with $T(0): = \min\{n\geq 1: X_n = 0\}$, the distribution of $T(0)$ with respect to  $Q_{\overline{\omega}}$ is stochastically dominated by the distribution of $T(0)$ with respect to $Q_ \omega$. 
We do this by coupling two Markov chains: the Markov chain $(\overline{X}_n)$ moves according to $Q_{\overline{\omega}}$, with $\overline{X}_0= 0$, the
Markov chain 
$(X_n)$ moves according to $Q_\omega$, with $X_0=0$, in such a way that $(\overline{X}_n)$ returns to $0$ before (or at the same time as) 
$(X_n)$. Let $(\overline{X}_n)$ move according to $Q_{\overline{\omega}}$.
If $\overline{X}_n \leq 0$ and $\overline{X}_{n+1} < \overline{X}_n$ and $X_n\leq 0$, then also $X_{n+1} < X_n$: this is possible since $Q_{\overline{\omega}}[\overline{X}_{n+1} < \overline{X}_n] = 1-M \leq Q_\omega[X_{n+1} <  X_n]$.
If $\overline{X}_n > 0$ and $\overline{X}_{n+1} > \overline{X}_n$ and $X_n>0$, then also $X_{n+1}> X_n$: this is possible since $Q_{\overline{\omega}}[\overline{X}_{n+1} > \overline{X}_n] = w \leq Q_\omega[X_{n+1} >  X_n]$. If $\overline{X}_n > 0$ and $\overline{X}_{n+1} > \overline{X}_n$ and $X_n< 0$, then $X_{n+1} <  X_n$: this is possible since $Q_{\overline{\omega}}[\overline{X}_{n+1} > \overline{X}_n] = w \leq 1-M \leq Q_\omega[X_{n+1} < X_n]$. 
Now, $(\overline{X}_n)$ visits $0$ before (or at the same time as) as $(X_n)$ does and this proves the statement. \\
Second step: we show that for each $n$, the distribution of $\xi(n,0)$ with respect to $Q_\omega$ is stochastically dominated by the distribution of 
$\xi(n,0)$ with respect to $Q_{\overline{\omega}}$. Let $T_1, T_2, \ldots $ be i.i.d. copies of $T(0)$. Then, $Q_\omega[\xi(n,0) \geq k] = P_\omega[\sum_{j=1}^k T_j \leq n] \leq Q_{\overline{\omega}}[\sum_{j=1}^k T_j \leq n] = Q_{\overline{\omega}}[\xi(n,0) \geq k]$, where the inequality follows from the first step.\\
Third step: For $x \in \integers$, the distribution of $\xi(n, x)$ with respect to $Q_\omega$ is dominated by the distribution of $\xi(n,0)$ with respect to  $Q_{\theta^x\omega}$ (more precisely, with $T(x): = \inf\{n: X_n = x\}$ denoting the first hitting time of $x$, the distribution of $\xi(n, x)$ with respect to $Q_\omega$ is the distribution of $\xi((n-T(x)), 0){\bf I}_{T(x) \leq n}$ with respect to  $Q_{\theta^x\omega}$). Now, apply the second step with $\theta^x\omega$ instead of $\omega$.\\
To show (\ref{nubar}), define, for $\omega \in [w, M]^\integers$, the corresponding potential $V$ as in Section \ref{s:intro}. Note that if $V$ is in the support of $\widetilde{V}$, the Markov chain defined by $Q_\omega$ is positive recurrent with invariant probability measure
\begin{equation}
\nu(x) =   \frac{\exp(- {V} (x-1))+ \exp( - {V} (x))}{ 2\sum\limits_{x \in \integers} \exp(- {V} (x))}\, ,
\end{equation}
and we have
$\nu(x) = \lim\limits_{n \to \infty}\frac{1}{n} \xi(n, x)$, hence (\ref{nubar}) follows from the stochastic domination.\qed

Intuitively, $\overline{V}$ is the steepest possible (infinite) valley which maximises the occupation time of $0$ (if $M \leq 1-w$).
Let $\varepsilon > 0$. Since $(X_n)$ is a positive recurrent Markov chain with respect to $Q_{\overline{\omega}}$, we have
$$
\lim_{n\to \infty} \, \frac{\xi(n,0)}{n} = \overline{\nu}(0)\qquad 
Q_{\overline{\omega}}\, \hbox{\rm -a.s.},  
$$
and
\begin{equation}
\label{expdec}
Q_{\overline{\omega}}\left[\frac{\xi(n,0)}{ n} \geq \overline{\nu}(0) + \varepsilon\right] \leq e^{-C(\varepsilon)n} ,
\end{equation}
where $C(\varepsilon)$ is a strictly positive constant depending only on $\varepsilon$. For $x \in \integers$, let $\theta_x$ be the shift on $\Omega$, i.e. $(\theta_x\omega)(y) = \omega(x+y)$, $y \in \integers$.
We have
\begin{eqnarray}
Q_\omega\left[ \frac{\xi^*(n)}{ n}\geq \overline{\nu}(0)+\varepsilon\right]
&\leq& \sum\limits_{x=-n}^n Q_{\theta^x\omega}\left[\frac{\xi(n,0)}{ n} \geq \overline{\nu}(0) + \varepsilon\right]\cr
&\leq& (2n+1) Q_{\overline{\omega}}\left[\frac{\xi(n,0)}{ n} \geq \overline{\nu}(0) + \varepsilon\right]\cr
&\leq& (2n+1) e^{-C(\varepsilon)n} ,
\end{eqnarray}
where we used Lemma \ref{coupling} for the second inequality and (\ref{expdec}) for the last inequality. By the Borel-Cantelli lemma, we conclude that
\begin{equation}
\label{upplimsup}
     \limsup_{n\to \infty} \, \frac{\xi^*(n)}{ n} \leq \overline{\nu}(0)
\qquad  Q_\omega \hbox{-a.s., for all  }  \omega \in [w, M]^\integers\, .
 \end{equation}
We have to show that (\ref{upplimsup}) holds true $\P$-a.s., i.e. for our RWRE on $\integers_+$, with reflection at $0$. In order to do so, we need the following modification of Lemma \ref{coupling} for environments on $\integers_+$.
Let $K \geq 1$ and let the environment $\overline{\omega}^{(K)}$ be defined as follows.
$$
{\overline \omega}^{(K)}_x : = \begin{cases}
 M, &  0< x \leq K \\
w, & x >  K\\
1, & x= 0. \end{cases}
$$
For the corresponding potential $\overline V^{(K)}$, $\exp(- {\overline V^{(K)}(x)})$ is given as follows: 
$$
\exp(- {\overline V^{(K)}(x)})  = \begin{cases}
\left(\frac{M}{1-M}\right)^x,  & 0< x \leq K \\
\left(\frac{M}{1-M}\right)^K\left(\frac{w}{1-w}\right)^{x-K},   & x > K \\
1, & x = 0. \\ 
\end{cases}
$$
The Markov chain on $\integers_+$ defined by $P_{\overline{\omega}^{(K)}}$ is positive recurrent with invariant probability measure ${\overline\nu}^{(K)}$ given by
\begin{equation}\label{nubarK}
{\overline\nu}^{(K)}(x):=
\begin{cases}
\frac{\exp(- {\overline V}^{(K)} (x-1))+ \exp( - {\overline V}^{(K)} (x))}{1+ 2\sum\limits_{x \geq 1} \exp(- {\overline V^{(K)}} (x))}, & x\geq 1\\
\frac{1}{1+ 2\sum\limits_{x \geq 1} \exp(- {\overline V^{(K)}} (x))}, & x=0\, .\end{cases}
\end{equation}

\begin{lemma}
\label{couplingrefl}
Assume $M \leq 1-w$. For all $K\geq 1$, if $\omega \in\{1\}\times [w, M]^{\integers_+}$ satisfies $\omega_0 =1$ and $b_n(\omega) \to \infty$, then for all $n \in \naturals$ such that $b_n(\omega) \geq 2K$, and all $x \in \integers$, $x\geq -b_n +K$, the distribution of $\xi(n,b_n+x)$ with respect to $P_\omega$ is stochastically dominated by the distribution of $\xi(n,K)$ with respect to 
$P_{\overline{\omega}^{(K)}}[\cdot\, |X_0 = K]$. Further, ${\overline\nu}^{(K)}(K) \to {\overline\nu}(0)$ for $K \to \infty$.
\end{lemma}
The proof of Lemma \ref{couplingrefl} is similar to the proof of Lemma \ref{coupling}. The last statement follows from (\ref{nubarK}).

Let $\varepsilon > 0$. Choose $K$ such that ${\overline\nu}^{(K)}(K) \leq  {\overline\nu}(0) + \frac{\varepsilon}{2}$.
Since $(X_n)$ is a positive recurrent Markov chain with respect to $P_{\overline{\omega}^{(K)}}[\cdot\,  |X_0 = K]$, we have
$$
\lim_{n\to \infty} \, \frac{\xi(n,K)}{n} = \overline{\nu}^{(K)}(K)\qquad 
P_{\overline{\omega}^{(K)}}\left[\cdot \Big| X_0 = K\right] \, \hbox{\rm - a.s.},  
$$
and
\begin{equation}
\label{expdecK}
P_{\overline{\omega}^{(K)}}\left[\frac{\xi(n,K)}{ n} \geq \overline{\nu}^{(K)}(K) + \varepsilon \Big| X_0  = K\right]  \leq e^{-C(\varepsilon)n} ,
\end{equation}
where $C(\varepsilon)$ is a strictly positive constant depending only on $\varepsilon$, and on $K$. For $P$-a.a. $\omega$, choose $n_0(\omega)$ such that for $n \geq n_0(\omega)$, $b_n(\omega) \geq 2K$ and $\frac{1}{n}\sum\limits_{x=0}^K \xi(n,x) \leq \frac{\varepsilon}{4}$, $P_\omega$-a.s. (this is possible since $b_n(\omega)$ increases to 
$\infty$, $P$-a.s. for $n \to \infty$, and $\frac{1}{n}\sum\limits_{x=0}^K \xi(n,x) \to 0$, $\P$-a.s. for $n \to \infty$ since $P_\omega$ is null-recurrent for $P$-a.a. $\omega$).
We then have $P$-a.s. for $n \geq n_0(\omega)$
\begin{eqnarray}
P_\omega\left[ \frac{\xi^*(n)}{ n}\geq \overline{\nu}(0)+\varepsilon\right]
&\leq& P_\omega\left[ \frac{\xi^*(n)}{ n}\geq \overline{\nu}^{(K)}(K)+\frac{\varepsilon}{2}\right]\cr
&\leq& \sum\limits_{x=-b_n+K}^n P_\omega\left[\frac{\xi(n,b_n+x)}{ n} \geq \overline{\nu}^{(K)}(K) + \frac{\varepsilon}{4}\right]\cr
&\leq& (2n+1) P_{\overline{\omega}^{(K)}}\left[\frac{\xi(n,K)}{ n} \geq \overline{\nu}^{(K)}(K) + 
\frac{\varepsilon}{4}\Big| X_0=K\right]\cr
&\leq& (2n+1) e^{-C(\varepsilon /4)n} ,
\end{eqnarray}
where we used Lemma \ref{coupling} for the third inequality and (\ref{expdecK}) for the last inequality. By the Borel-Cantelli lemma, we conclude that
\begin{equation}
\label{upplimsuprefl}
     \limsup_{n\to \infty} \, \frac{\xi^*(n)}{ n} \leq \overline{\nu}(0)
\qquad \hbox{\rm $P_\omega$ -a.s., for  } P \hbox{- a.a. } \omega 
 \end{equation}

Together with (\ref{nubar}), this finishes the proof in the case $M \leq 1-w$: the value of $c$ is computed easily from (\ref{nubar}).
In the case $M > 1-w$, Lemma \ref{coupling} is replaced with the following
\begin{lemma}
\label{coupling2}
Let the environment $\overline{\omega}$ be defined as above. Assume $M > 1-w$. For {\sl all} $\omega \in [w, M]^\integers$, all $x \in \integers$ and all $n \in \naturals$, the distribution of $\xi(T(1)+n,x)-\xi(T(1),x) $ with respect to $Q_\omega$ is stochastically dominated by the distribution of $\xi(T(1)+n,1)$ with respect to $Q_{\overline{\omega}}$. In particular, 
\begin{equation}\label{nubar2}
\sup \{z: z \in \hbox{\rm supp} \, \sup_{x\in \integers} {\nu}(x)\}
= {\overline\nu}(1)=
\frac{\exp(- {\overline V} (0))+ \exp( - {\overline V} (1))}{ 2\sum\limits_{x \in \integers} \exp(- {\overline V} (x))} . 
\end{equation}
\end{lemma}
The rest of the proof is analogous to the case $M \leq 1-w$.
\qed

\section{Further remarks}
\label{s:remarks}

1. Consider RWRE on the integers axis, i.e. with transition probabilities
$Q_\omega [X_{n+1} = x+1 \, |\, X_n =x] = \omega_x = 1-Q_\omega [X_{n+1} = x-1 \, | \, X_n =x]$ for $x \in \integers $.
Then,
Theorem \ref{empdistr} can be modified in the following way (we refer to \cite{bil} for details). First, we have to replace $b_n$ in (\ref{mainconv}) with $\widehat b_n$ defined as follows.
We call a triple $(a,b,c)$ with $a < b < c$ a {\sl valley} of $V$ if
$$
V(b) = \min\limits_{a \leq x \leq c }V(x),\quad V(a) = \max\limits_{a \leq x \leq b }V(x),\quad 
V(c) = \max\limits_{b \leq x \leq c }V(x)\, .
$$
\noindent The depth of the valley is defined as $d_{(a,b,c)} = (V(a) - V(b)) \wedge (V(c) - V(b))$.
Call a valley $(a, b, c)$ {\sl minimal} if for all valleys $(\widetilde{a},\widetilde{b},
\widetilde{c})$ with $a < \widetilde{a}$, $\widetilde{c} < c$, we have 
$d_{(\widetilde{a},\widetilde{b},\widetilde{c})} < d_{(a,b,c)}$. Consider the smallest minimal valley $(\widehat a_n, \widehat b_n, \widehat c_n)$ with $\widehat a_n < 0 < \widehat c_n$ and $d_{(a_n,b_n,c_n)} \geq \log n+ (\log n)^{1/2}$, i.e. $(\widehat a_n,\widehat b_n, \widehat c_n)$ is a minimal valley with $d_{(\widehat a_n,\widehat b_n, \widehat c_n)} \geq \log n+ (\log n)^{1/2}$ and for every other valley $(\widetilde a_n, \widetilde b_n, \widetilde c_n)$ with these properties, we have $\widetilde a_n < \widehat a_n$ or $\widehat c_n < \widetilde c_n$.
If there are several such valleys, take
$(\widehat a_n, \widehat b_n ,\widehat c_n)$ such that $|\widehat b_n|$ is minimal (and $\widehat b_n >0$ if there are two possibilities). 
Further, one has to replace the limit measure $\nu$ in  (\ref{mainconv}) with $\widehat \nu$ defined in as follows.
Let $\widetilde V_{\rm left}= (\widetilde V_{\rm left}(x), \, x\in \integers)$ be a collection of random variables distributed as $V$ conditioned to stay strictly positive for $x > 0$ and non-negative for $x<0$.
(Recall that ${\widetilde V}= ({\widetilde V}(x), \, x\in \integers)$ is a collection of random variables distributed as $V$ conditioned to stay non-negative for $x > 0$ and strictly positive for $x<0$.)
Set
\begin{equation}\label{nuleftdef}
    \nu_{\rm left}(x)
:= \frac{\exp(- \widetilde V_{\rm left} (x-1))+ \exp( -
    \widetilde V_{\rm left}(x))}{ 
    2\sum\limits_{x \in \integers} \exp(-
    \widetilde V_{\rm left}(x))}\, ,\qquad x \in \integers.
\end{equation}
In other words, $\nu_{\rm left}$ is defined in the same way as $\nu$, replacing $\widetilde V$ in (\ref{nudef}) with $\widetilde V_{\rm left}$.
Now, let $Q$ be the distribution of $\nu$ and $Q_{\rm left}$ be the distribution of $\nu_{\rm left}$ and let $\widehat \nu$ be a random probability measure on $\integers$ with distribution  $\frac{1}{2} \left( Q + Q_{\rm left}\right)$.
Then, (\ref{mainconv}), (\ref{fucoll}), (\ref{sil}) and (\ref{maxloc}) hold true for RWRE on the integer axis, if we replace $b_n$ with $\widehat b_n$ and $\nu$ with $\widehat\nu$. \\
Since the right hand side of (\ref{theconst}) is a function of $M$ and $w$, say $\varphi(M, w)$, for which we have $\varphi(M, w) = \varphi(1-w, 1-M)$, we see that replacing $V$ with $V_{\rm left}$ in (\ref{stardef}) yields the same constant $c$ as in (\ref{cocon}) and therefore Theorem \ref{t:LIL} carries over verbatim.\\
2. We showed that $\limsup_{n\to \infty} \frac{\xi^*(n)}{ n}$ is $\P$-a.s. a strictly positive constant and gave its value.
We have
 $\liminf_{n\to \infty} \frac{\xi^*(n)}{ n} =0$ $\P$-a.s. which is a consequence of (\ref{maxloc}). It was shown in \cite{ANYZ} that there is a strictly positive constant $a$ such that
$\liminf_{n\to \infty} \frac{\xi^*(n)\log\log\log n}{ n} = a$, $\P$-a.s.
The value of $a$ is not known.

{\bf Acknowledgement} We thank Gabor Pete, Achim Klenke and Pierre Andreoletti for pointing out mistakes in earlier versions of 
this paper and Amir Dembo for discussions. We also thank the referee for several valuable comments.

\end{document}